\documentclass[preprint,12pt]{elsarticle}

\usepackage{graphics}
\usepackage{graphicx}
\usepackage{epsfig}

\journal{???}

\begin{document}
\begin{frontmatter}

\title{House of Graphs: a database of interesting graphs}


\author[gent]{Gunnar Brinkmann}
\ead{Gunnar.Brinkmann@UGent.be}
\author[gent]{Kris Coolsaet}
\ead{Kris.Coolsaet@UGent.be}
\author[gent]{Jan Goedgebeur}
\ead{Jan.Goedgebeur@UGent.be}
\author[mons] {Hadrien M\'elot}
\ead{Hadrien.Melot@UMons.ac.be}

\address[gent]{Applied Mathematics \& Computer Science\\
  Ghent University\\
Krijgslaan 281-S9, \\9000 Ghent, Belgium\\ }
\address[mons]{Institut d'Informatique\\
Universit\'e de Mons\\
Place du Parc, 20\\
7000 Mons, Belgium }

\begin{abstract}
In this note we present {\em House of Graphs} (\verb+http://hog.grinvin.org+)
which is a new database of graphs. 
The key principle is to have a searchable database 
and offer -- next to complete lists of some graph classes --
also a list of 
special graphs that already turned out to be {\em interesting}
and {\em relevant} in the study of graph theoretic problems
or as counterexamples to conjectures.
This list can be extended by users of the database.
\end{abstract}


\end{frontmatter}


\section{Introduction}

On the web one can find several websites with lists of graphs 
-- see e.g.\ 
the lists of
Brendan McKay \cite{graphlists_brendan}, Gordon Royle \cite{graphlists_gordon},
Markus Meringer \cite{graphlists_markus}, Frank Ruskey \cite{graphlists_frankruskey}
or Ted Spence \cite{graphlists_tedspence}.
In the {\em Atlas of Graphs} \cite{atlas_of_graphs}
one can even find hardcopy lists with lots of pictures of graphs. 
Graphs are given as adjacency lists, sometimes with, sometimes without additional data
and in various formats.
Such lists can serve as a 
source for intuition when one studies some conjecture and even as a possible source for counterexamples.
The website ChemSpider \cite{chemspider} specialising in chemical structures 
not only gives lists of (graph) structures but also allows
several chemically motivated ways to search the database.

As the number of graphs grows very fast (see sequence A001349 in the
On-line Encyclopedia of Integer Sequences \cite{OEIS}), 
 even with modern
computer systems one will not be able to test conjectures or study
properties on {\bf all} graphs up to $n$ vertices unless $n$ is very small.
For example, the complete list of connected graphs with 14 vertices contains
already 29.003.487.462.848.061 graphs -- so that even the computation of
simple invariants for all of these graphs is quite a challenge.
If one restricts the class of graphs one wants to study or test (e.g., only regular
graphs, only 3-regular graphs, etc.), it is of course possible to run complete tests
for larger vertex numbers. Unfortunately this is not always possible and 
most restricted classes of graphs also
grow exponentially fast imposing a relatively small limit on the
number of vertices for which complete lists can still be tested.
The main new functionality of the database presented here is to offer a list
of relatively small size that still gives you a good chance to find
counterexamples or obtain results with your tests that allow to judge
the general situation. For these lists there is also some additional functionality
that will be described later on.

Some graphs (e.g., the famous
Petersen graph or the Heawood (3,6)-cage on 14 vertices) or graph classes (e.g., snarks)
appear repeatedly in the literature while others will probably always be just part
of the huge mass. So similar to Orwell's famous words: all graphs are interesting, but 
some graphs are more interesting than others.

If one wants to test a conjecture on a list of graphs, the ideal case would be if one 
could restrict the tests to graphs that are {\em interesting} or {\em relevant}
for this conjecture. Here the meaning of {\em interesting} and {\em relevant}
is very vague, but this already shows how much the question whether a graph
is interesting or not depends on the question one wants to study.

In this paper we will not try to give an exact definition of the terms {\em interesting}
or {\em relevant} -- not even in the restricted form {\em interesting/relevant for
a certain invariant}. If you think of a specific graph, this is sufficient to separate
the graph from the huge mass of other graphs, proving that the graph is interesting 
in some respect -- and of course the database allows you to add that graph and also
offers the possibility to say for which invariants the graph is especially interesting.
We consider e.g.\ counterexamples to known conjectures 
and extremal graphs to be interesting.

In order to be a rich source of possible counterexamples right from the start,
we already added 1570 graphs to the database. Most of these graphs are extremal graphs
found by {\em GraPHedron} \cite{graphedron}, but there are
also other graphs like maximal
triangle free graphs that occured as Ramsey graphs for triangle Ramsey numbers \cite{BrBrHa96_2}
or graphs that occured as named graphs in the mathematical literature. 

The graphs found by {\em GraPHedron} are extremal in the following sense:
for a given set of $p$
graph invariants, GraPHedron uses a polyhedral approach placing the 
graphs in $p$-dimensional space and computing the facets that bound the set.
These facets determine linear inequalities between the graph invariants
for the set of graphs under consideration. Extremal graphs are graphs where
the corresponding point in $p$-space is a vertex of the polyhedron. This means
that the extremal graphs determine the inequalities and in the database they
are listed as interesting with respect to the invariants occurring in the 
inequality.

In some cases a lot of graphs are situated at the same extremal point.
We call these sets of graphs {\em conglomerates}. 
Typical conglomerates are e.g.\ trees or complete
bipartite graphs. Adding all conglomerates would increase the number of
graphs dramatically and would lead to lists that are much larger than desired.
All graphs in conglomerates share a specific structure that made them extremal
points for the given set of invariants. But this specific structure can be represented
by any element of a conglomerate.
So we decided to represent conglomerates by one of their elements.
Choosing a minimal set of graphs representing all
conglomerates is in fact a Minimum Set Covering problem and therefore NP-complete,
so we have used some greedy heuristic to determine a relatively small set 
of graphs representing all conglomerates that occurred during the 
computations.

\section{Functionality of the website}

A basic functionality one must expect from every database of graphs
is the possibility to download some specific lists of graphs. At the moment 
the website {\em House of Graphs}
offers (among others) the following lists:

\begin{itemize}

\item All graphs registered as interesting in the database.

\item All snarks up to 34 vertices (girth 4) resp. 36 vertices (girth 5).

\item All IPR-fullerenes up to 160 vertices.

\item Complete lists of regular graphs for various combinations of
degree, vertex number and girth.

\item Vertex-transitive graphs.

\item Some classes of planar graphs.

\end{itemize}

Some of these lists are physically situated on the same server as 
the website itself, but others are just links to other peoples websites, 
like those of
Brendan McKay \cite{graphlists_brendan}, Gordon Royle \cite{graphlists_gordon},
Ted Spence \cite{graphlists_tedspence}
or Markus Meringer \cite{graphlists_markus}. 

In some cases it is much faster
to generate the graphs instead of storing the list and downloading it.  
{\em House of Graphs}
offers source codes
of some graph generation programs (like e.g.\ {\em snarkhunter} \cite{tricycle},
{\em minibaum} \cite{Br96_1} or {\em MTF} \cite{BrBrHa96})
and contains links to other pages where
such programs (like {\em genreg} \cite{MM98}, {\em geng }\cite{McK96} or {\em FreeTree } \cite{freetree})
can be downloaded.

\subsection{Search functionality of the website}

The list of graphs that are marked as {\em interesting} plays a special role 
in this database. For graphs in this list, a lot of invariants 
(like the chromatic number, the clique number, the diameter, etc.) 
and also embeddings (drawings)
are precomputed and stored. 

Furthermore this database can be searched
in various ways. In this context ``searching'' always means that the
list of all graphs in the database is reduced to the sublist of graphs that satisfy
the search criterion.

\medskip

It can be searched for:

\begin{itemize}

\item 
A keyword -- e.g.\  a search for 
``Petersen'' will give a list containing the Petersen graph, but also a graph that is a
Kronecker product involving the line graph of the Petersen graph as that involvement 
is mentioned in the text belonging to the graph.

\item Graphs
that are marked as interesting for a certain invariant (e.g.\  all graphs that are 
marked interesting for the matching number).

\item Graphs where 
the value of an invariant is in a certain range (e.g.\  graphs for which the matching number
$m(G)$ is between 11 and 13).

\item A specific graph to see whether it is already in the database
and what information is given on it.
The graph can either be
described by uploading a file in {\em graph6} format or an easy text format that is described 
on the site. Small graphs can also be described by drawing them using a simple
graph drawing applet on the site.

\end{itemize}

Future plans are to make it also possible to restrict the list to graphs that satisfy
or do not satisfy some expression in the given invariants -- e.g.\  to graphs that fulfil  
$m(G) \le (|V(G)|/2) - 2$.

By applying several restriction steps, these criteria can be combined with each other.

For example, to search for the cubic 6-cage one could apply the following restriction steps: girth 6 $\rightarrow$ regular $\rightarrow$ average degree 3 $\rightarrow$ 14 vertices.

\begin{figure}[h!t]
	\centering
	\includegraphics[width=0.4\textwidth]{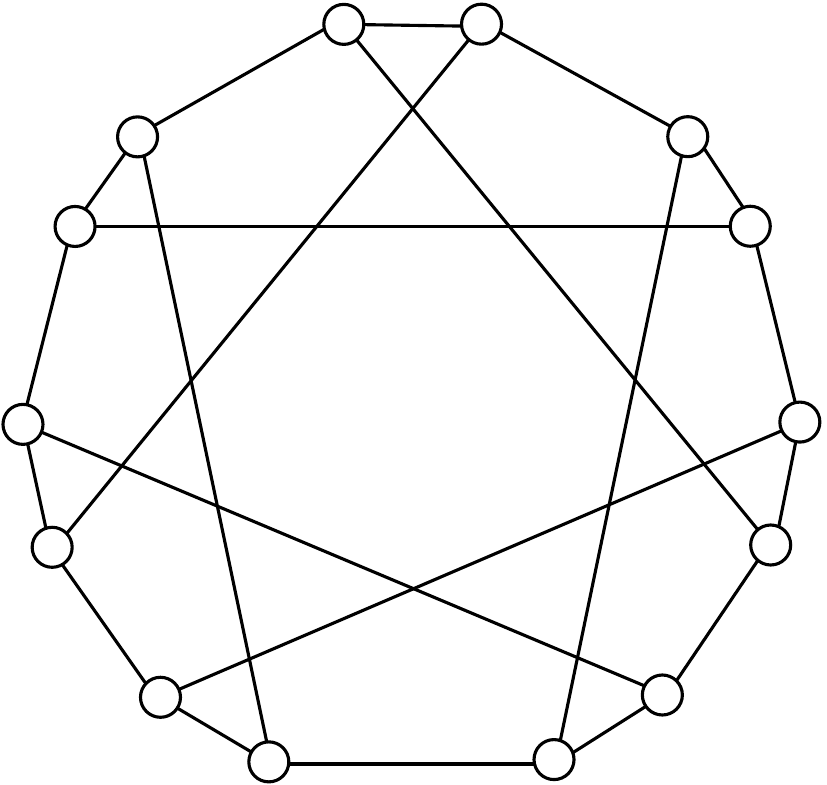}
	\caption{The Heawood graph.}
	\label{fig:heawood_graph}
\end{figure}

This results in the graph in Figure~\ref{fig:heawood_graph}, which is also known as the Heawood graph. The same graph is found by searching for the keywords ``Heawood'' or ``(3,6)-cage''.

All lists -- complete lists or lists 
that are the result of a search --
can be downloaded in {\em graph6} format or the binary multicode format which are explained on the
website. They can also be downloaded in a self-explanatory human readable text format.
By browsing through the lists and clicking on graphs in the list additional information
on the graph is displayed -- like invariant values, the name of the graph, for which invariants
the graph is considered interesting or comments on the graph.

\subsection{Submitting graphs to the site}

It is essential that the list of graphs marked as {\em interesting} can grow whenever 
graphs are discovered that turn out to be relevant in the study of some graph theoretic
problem. 
To this end it is necessary that people 
can submit new graphs to the list.

In order to submit a graph to the site you have to register first (which is also free). 
Having registered,
a graph can be submitted in several ways: it can be submitted in any of the
formats that are also available as download formats for the graphs. Furthermore some 
basic graph editor allows to draw graphs. 
After the graph has been inputted, it is checked whether the graph
(or an isomorphic copy) is already in the list. 

If the graph is new, the user submitting the graph is considered the owner of
the graph and some information on the
graph can be filled in, like e.g.\ the name of the graph, 
whether it is a counterexample 
to a conjecture, whether it is the result of a construction, or references to publications 
where it first occurred. Furthermore, the user can mark
some invariants for which he thinks the graph must be considered interesting.  
In case the graph was submitted by drawing it, this drawing is
used as the default visual presentation of the graph, otherwise a simple automatic
drawing is computed.  

If the graph
was already in the list, it is still possible to add comments to the existing graph
e.g.\ about properties not mentioned by the owner when the graph was submitted.
The basic information about the graph -- e.g.\ the name of the graph -- can only be
changed by the owner.

After the graph has been submitted, processes to compute the invariant values
are submitted to a background queue. The routines to compute the invariants are from
the graph theory environment {\em Grinvin} \cite{grinvin}.
Depending on the load of the server
and the difficulty of the invariant to be computed,
the time until the invariant information is displayed can range from 
milliseconds to days. For very large graphs and NP-complete invariants the invariant value
may even remain unknown. 
Information on known values of invariants for graphs where the computation of
the invariants may take very long, can be given as a comment.

In the future it should also be possible to submit graphs by uploading
e.g.\ pdf or jpg files of the drawing of the graph. At the moment this is not possible.

\section*{Acknowledgements}

Jan Goedgebeur is supported by a PhD grant from the Research Foundation of Flanders (FWO).


\end{document}